\documentclass[a4paper,12pt]{scrartcl}

\usepackage[utf8]{inputenc}

\usepackage{amsmath,amsthm,graphicx}
\usepackage[commonsets]{skmath}
\usepackage{bm}

\usepackage{graphicx,color}

\usepackage{fnpct}

\usepackage{bera}
\usepackage[charter]{mathdesign}

\usepackage{cite}
\usepackage{xspace,microtype}

\begin{document}

\title{A short note on the stability of a class of parallel
  Runge-Kutta methods}

\author{Friedemann Kemm}

\maketitle

\begin{abstract}
  With this short note, we close a gap in the linear stability theory
  of block predictor-corrector Runge-Kutta schemes originally proposed
  for the parallel solution of ODEs. 
\end{abstract}

\section{Introduction}
\label{sec:introduction}

The advent of parallel computers, especially of machines with shared
memory, gave rise to the evolution of block parallel Runge-Kutta
methods, especially of predictor-corrector type,
cf.~\cite{HouwSommparit,pirkas,voss-abbas-1997,burrage-suhartanto-97,rosser}
to name just a few studies in that field. A class that caught some
attention starts with explicit Euler as a predictor which is followed
by several corrector steps, all with the same method. The according
Butcher Tableau is of the form
\begin{equation}
  \label{eq:10}
  \begin{array}{c|ccccc}
    \vec 0 & & & & & \\
    \vec c & \vec A & & & & \\
    \vec c & \vec 0 & \ddots & & & \\
    \vdots & & & \ddots & & \\
    \vec c & \vec 0 & \dots & \vec 0 & \vec A & \\
    \hline\\[-.6em]
           & \vec 0^T & \dots & \dots & \vec 0^T & \vec b^T
  \end{array}
\end{equation}
where we used the same symbol for the zero matrix and the zero vector
of dimension~\(s\). The matrix~\(\vec A\) as well as the
vectors~\(\vec b\) and~\(\vec c\) are taken from the representation of
the corrector method
\begin{equation}\label{eq:11}
  \begin{array}{c|c}
    \vec c &  \vec A \\
    \hline\\[-.6em]
          & \vec b^T
  \end{array}\quad\text{with}\quad
  \vec A \in \mathbb R^{s\times s},\ \vec b, \vec c \in \mathbb R^s\;,\
  c_i\in[0,1]\ \forall i\;.
\end{equation}
For this type of method in each corrector step~\(s\)~fluxes can be
computed simultaneously without any coupling. Thus, these methods are
well suited to parallel computers with shared
memory. In~\cite{HouwSommparit}, Houwen and Sommeijer discuss the
stability up to the optimal order---when the order of the corrector
method is achieved with minimal number of correction steps---and a
step-size control that makes use of the nature of the
corrector method.
In the context of PDEs it is unusual to rely on these control
strategies. The usual way is to determine the size of the time steps
by the linear stability of the method. As a consequence, the question
arises if more correction steps than needed for the optimal order
would give a benefit via an increased region of absolute linear
stability.

In this study, we close a gap in the linear stability theory of
methods of type~\eqref{eq:10}. We show that the stability function is
always completely determined by the stability function of the corrector
method. It is always a Taylor polynomial of the corrector's stability
function.

The paper is organized as follows: In the next section, we prove the
main result. After this, we discuss some consequences for the
application of these methods, especially in the context of PDEs. The
study is then concluded with a short summary and outlook. 
For the general theory of Runge-Kutta methods, we refer to the standard
textbooks~\cite{HW1,HW2,butcherbuch}. 

\section{The main result}
\label{sec:main-result}

In this section, we prove the main result that the linear stability
function~\(R_{\text{PC},\,m}\), where~\(m\) denotes the number of
corrections,  of method~\eqref{eq:10} is a truncated Taylor series
expansion (with respect to the origin) of the stability function~\(R\) of
the corrector~\eqref{eq:11}. 

In the following, the notation~\(\vec 1 = (1,\dots,1)^T\) is used
without reference to the actual dimension as is also done with the
identity matrix~\(\vec I\). With this notation, the stability function
of the corrector can be written as
\begin{equation}
  \label{eq:5}
  R(z) = 1 + z\vec b^T (\vec I - z \vec A)^{-1} \vec 1\;. 
\end{equation}


For method~\eqref{eq:10}, Houwen and Sommeijer~\cite{HouwSommparit} computed the
linear stability function, which in the case of~\(m\)~correction steps
is\footnote{They originally restricted their paper to correctors where
  the last stage value is already the new value. But this restriction
  is only needed for their step-size control.} 
\begin{equation}
  \label{eq:6}
  R_{\text{PC},\,m}(z) = 1 + z\vec b^T\vec 1 + z^2\vec b^T \vec A \vec 1 +
  \dots + z^{m+1}\vec b^T \vec A^m \vec 1\;. 
\end{equation}
They proof that for~\(m=p-1\), this is just the according Taylor
polynomial of the exponential function
\begin{equation}
  \label{eq:7}
  R_{\text{PC},\,p-1}(z) = 1 + z + \frac{1}{2!}z^2 + \frac{1}{3!}z^3 + \dots +
  \frac{1}{p!}z^p\;.   
\end{equation}
Their proof also implies that for~\(m<p-1\) the stability
function is a Taylor polynomial of degree~\(m+1\). Now we want to prove
that this is also the case for~\(m\geq p\).

If in general a matrix~\(\vec M\) depends on a parameter~\(z\),
i.\,e.~\(\vec M = \vec M(z)\), and is differentiable with respect
to~\(z\) in some open set~\(\Omega \subseteq \C\) and in addition
invertible for all~\(z\in\Omega\) then we have in~\(\Omega\)
\begin{equation}
  \label{eq:1}
  \pd{\vec M^{-1}}{z} = -\vec M^{-1} \pd{\vec M}{z} \vec M^{-1}\;. 
\end{equation}
Now we consider~\(\vec M = \vec I - z \vec A\) with a constant matrix~\(\vec A\) as above.
Apparently, in this case,~\(\vec M\) is differentiable in~\(\C\) and
invertible in some open neighbourhood~\(U\) of~\(0\). Thus, we have in~\(U\)
\begin{align}
  \label{eq:2}
  \pd{ (\vec I - z \vec A)^{-1}}{z} & = - (\vec I - z \vec
  A)^{-1} \pd{(\vec I - z \vec A)}{z}  (\vec I - z \vec A)^{-1}\\
  & =  (\vec I - z \vec A)^{-1} \vec A\, (\vec I - z \vec A)^{-1}\;. 
\end{align}
For higher derivatives, this can now be used recursively in connection
with the product rule. E.\,g.\ the second derivative can be computed as
\begin{align}
  \pd{ (\vec I - z \vec A)^{-1}}{z^2}
  & =  \pd{ (\vec I - z \vec A)^{-1}}{z} \vec A\, (\vec I - z \vec
    A)^{-1} + (\vec I - z \vec A)^{-1} \vec A\, \pd{ (\vec I - z \vec
    A)^{-1}}{z} \\
  & = 2 \cdot  (\vec I - z \vec A)^{-1} \vec A\, (\vec I - z \vec
    A)^{-1} \vec A\, (\vec I - z \vec A)^{-1}\;,
\end{align}
which recursively leads to
\begin{equation}
  \label{eq:3}
  \pd{ (\vec I - z \vec A)^{-1}}{z^n}
  = n! \cdot (\vec I - z \vec A)^{-1} \prod_{j=2}^n \vec A\, (\vec I
  - z \vec A)^{-1}\;. 
\end{equation}
Evaluated in~\(z=0\), we have
\begin{equation}
  \label{eq:4}
  \pd{}{z^n}R(0) = n! \vec b^T \vec A^{n-1}\vec 1\;,\qquad n=1,2,\dots 
\end{equation}

Together with equation~\eqref{eq:6} we get
\begin{equation}
  \label{eq:8}
  \begin{split}
    T_{R,\,m}(z) & = \sum_{n=0}^m \pd{}{z^n}R(0)\,z^n\\ 
    & = \sum_{n=0}^m \frac{1}{n!}n! \vec b^T \vec
    A^{n-1}\vec 1\,z^n\\ 
    & = 1 + z\vec b^T\vec 1 + z^2\vec b^T \vec A \vec 1 +
    \dots + z^{m+1}\vec b^T \vec A^m \vec 1\\
    & = R_{\text{PC},\,m}(z)\;. 
  \end{split}
\end{equation}
This means that in general, the stability function resulting
from~\(m-1\) corrector steps is the Taylor polynomial of degree~\(m\)
of the stability function of the method used as the corrector.

\section{Consequences for the application of these predictor-corrector
  methods}
\label{sec:cons-appl-these}

The result of the previous section has interesting consequences for
the practical use of these block predictor-corrector Runge-Kutta
schemes. It is possible to get the stability function and, thus, the
region of absolute stability directly from the stability function of
the corrector method, not only for the optimal order. In some cases,
additional corrector steps might further increase the stability
region and thereby the possible step sizes in a way that the
computational costs decrease. Note that in
the simulation of PDEs and systems of PDEs, e.\,g.\ in computational
fluid dynamics, the time steps are not chosen by accuracy, but by
stability. Furthermore, the order of the time integration method does
not need to exceed the order of the space discretization. 

As an example, we consider the 2-stage Radau~IIA method
\begin{equation*}
  \begin{array}{c|cc}
    1/3 & 5/12 & -1/12\\
    1   & 3/4  &  1/4 \\
    \hline\\[-.6em]
        & 3/4 & 1/4 \\
  \end{array}
\end{equation*}
which has the stability function
\begin{equation}\label{eq:9}
  R(z) = {\frac {1+z/3}{1-2/3\,z+1/6\,{z}^{2}}}\;.
\end{equation}
The poles of this function are
\begin{equation*}
  z_{1/2} = 2 \pm \ii\sqrt{2}\;. 
\end{equation*}
Both poles are sufficiently far from the stability region of the
optimal order method. We expect the Taylor series expansion of~\(R\)
with respect to~\(z_0 = 0\) to converge to~\(R\) in an open disc of
radius~\(\sqrt 6\) centered at the origin. Since Radau~IIA is
L-stable, we can expect numbers~\(m>3\) such that the
predictor-corrector method with~\(m-1\) corrections has a larger
region of stability than with~\(m=3\), which would already give us the
optimal order.

\begin{figure}
  \centering
  \includegraphics[width=.5\linewidth]{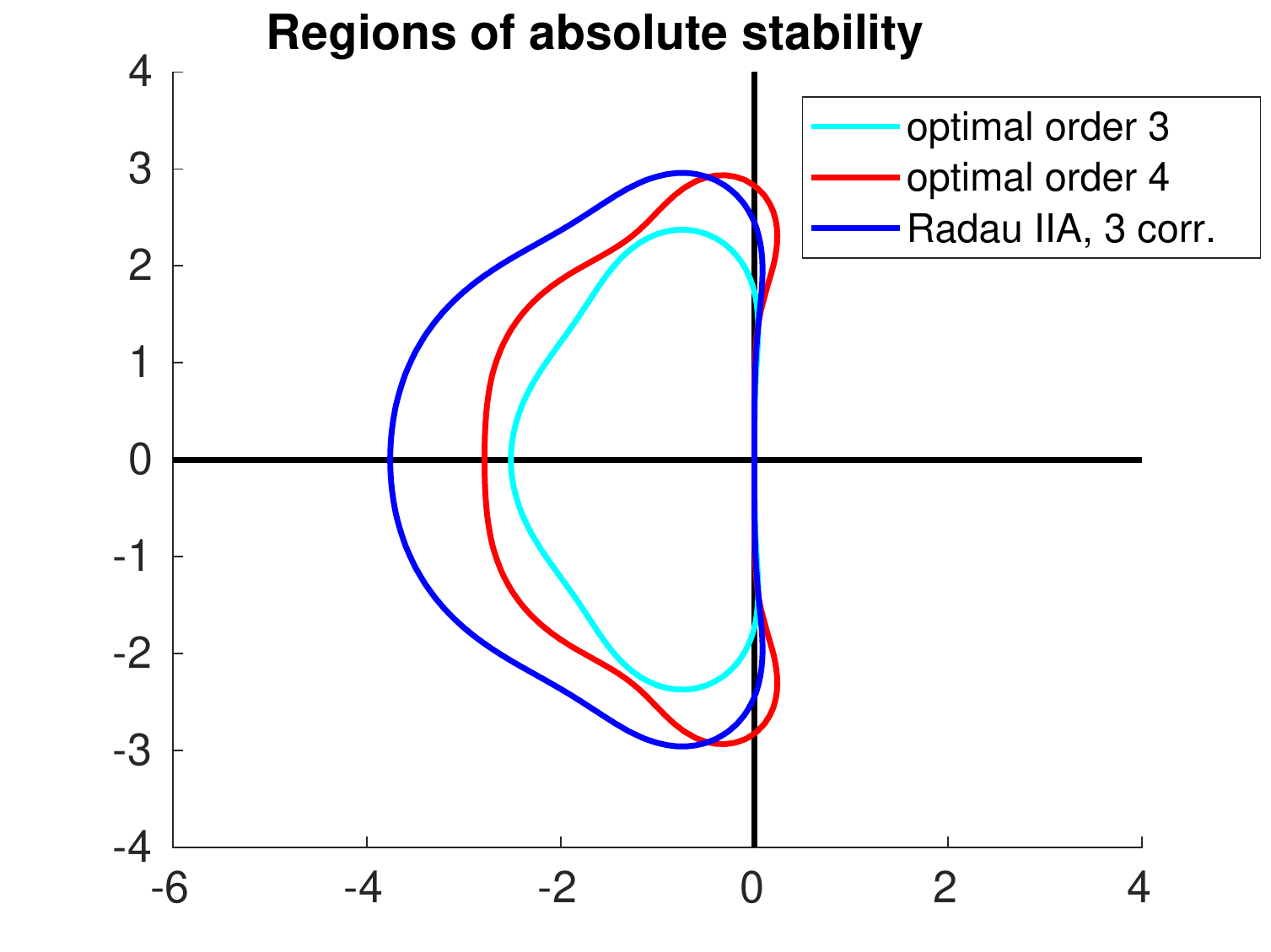}
  \caption{Regions of absolute stability for optimal orders~3 and~4
    and for the method with 2-stage Radau~IIA as corrector and
    3~corrector steps.}
  \label{fig:vglstab}
\end{figure}

In Figure~\ref{fig:vglstab}, we display the stability regions for the
optimal order with 2-stage Radau~IIA as a corrector, for optimal order
with a fourth order 2-stage corrector (Hammer-Hollingsworth), and with
Radau~IIA and an additional corrector step. it is obvious that, except
when order four is desired, it is advantageous to employ Radau~IIA
with three corrector steps. The number of (parallel) flux evaluations
per time span can be decreased.

By our main result, we know that this would also be true if we
replaced Radau~IIA by Radau~IA since both methods have the same
stability function.

\section{Conclusions and outlook}
\label{sec:conclusions-outlook}

With this study, we closed a gap in the linear stability theory of
block parallel predictor-corrector methods. A still open question is
what would happen if we replace the predictor, in our case explicit
Euler, by an explicit method with higher order. Will the stability
function again be uniquely determined by the stability function of the
corrector? Will it again yield truncated Taylor series expansions of the
correctors stability function? If so, our considerations regarding the
application to PDEs could also be transferred.


\bibliographystyle{amsplain} \bibliography{rkpcstab}

\providecommand{\bysame}{\leavevmode\hbox to3em{\hrulefill}\thinspace}
\providecommand{\MR}{\relax\ifhmode\unskip\space\fi MR }
\providecommand{\MRhref}[2]{%
  \href{http://www.ams.org/mathscinet-getitem?mr=#1}{#2}
}
\providecommand{\href}[2]{#2}
\begin{thebibliography}{1}

\bibitem{burrage-suhartanto-97}
K.~Burrage and H.~Suhartanto, \emph{Parallel iterated methods based on
  multistep {Runge}-{Kutta} methods of {Radau} type}, Adv. Comput. Math.
  \textbf{7} (1997), no.~1-2, 37--57 (English).

\bibitem{butcherbuch}
J.~C. Butcher, \emph{Numerical methods for ordinary differential equations.},
  2nd revised ed. ed., Hoboken, NJ: John Wiley \& Sons, 2008.

\bibitem{HW1}
Ernst Hairer, Syvert~P. N{\o}rsett, and Gerhard Wanner, \emph{Solving ordinary
  differential equations. {I}: {Nonstiff} problems.}, 2. rev. ed. ed., Springer
  Ser. Comput. Math., vol.~8, Berlin: Springer-Verlag, 1993.

\bibitem{HW2}
Ernst Hairer and Gerhard Wanner, \emph{Solving ordinary differential equations.
  {II}: {Stiff} and differential-algebraic problems.}, 2nd rev. ed. ed.,
  Springer Ser. Comput. Math., vol.~14, Berlin: Springer, 1996.

\bibitem{rosser}
J.~Barkley Rosser, \emph{A {Runge}-{Kutta} for all seasons}, SIAM Rev.
  \textbf{9} (1967), 417--452 (English).

\bibitem{HouwSommparit}
P.~J. van~der Houwen and B.~P. Sommeijer, \emph{Parallel iteration of
  high-order {Runge}-{Kutta} methods with stepsize control}, J. Comput. Appl.
  Math. \textbf{29} (1990), no.~1, 111--127 (English).

\bibitem{pirkas}
P.~J. van~der Houwen, B.~P. Sommeijer, and W.~A. van~der Veen, \emph{Parallel
  iteration across the steps of high-order {Runge}-{Kutta} methods for nonstiff
  initial value problems}, J. Comput. Appl. Math. \textbf{60} (1995), no.~3,
  309--329 (English).

\bibitem{voss-abbas-1997}
D.~Voss and S.~Abbas, \emph{Block predictor-corrector schemes for the parallel
  solution of {ODEs}}, Comput. Math. Appl. \textbf{33} (1997), no.~6, 65--72
  (English).

\end{thebibliography}

\end{document}